\documentclass[review,3p,times]{elsarticle}

\usepackage[russian, english]{babel}
\usepackage{amsthm,amssymb,amsmath,amsfonts}
\usepackage{graphicx}
\usepackage{hhline,array}
\usepackage{color}
\numberwithin{equation}{section} \numberwithin{figure}{section}

\newtheorem{theorem}{\bf Theorem}[section]

\usepackage{ecrc}

\volume{00}

\firstpage{1}

\runauth{I. R. Muftahov. A. N. Tynda, D. N. Sidorov}

\jid{APNUM}

\usepackage{amssymb}
\usepackage{amsthm}

\usepackage[figuresright]{rotating}

\begin{document}

\begin{frontmatter}



\dochead{}

\title{Numerical solution of Volterra integral equations of the first kind with discontinuous kernels}

\author[label1]{Ildar Muftahov}
\author[label2]{Aleksandr Tynda}
\author[label1,label3]{Denis Sidorov}

\address[label1]{Irkutsk National Research Technical University}
\address[label2]{Penza State University}
\address[label3]{Energy Systems Institute of Russian Academy of Sciences}

\begin{abstract}
{We propose the numerical methods for solution of the weakly regular linear and nonlinear evolutionary (Volterra) integral equation of the
first kind. The kernels of such equations have jump discontinuities along the
continuous curves (endogenous delays) which starts at the origin. In order to linearize these equations we use the modified Newton-Kantorovich iterative process.
Then for linear equations we propose two direct quadrature methods based on the piecewise constant and piecewise linear approximation of the exact solution.
The accuracy of proposed numerical methods is $\mathcal{O}(1/N)$ and $\mathcal{O}(1/N^2)$ respectively.
We also suggest a certain iterative numerical scheme enjoying the regularization properties. Furthermore, we adduce generalized numerical method for nonlinear equations.
  We employ the midpoint quadrature rule in all the cases.
In conclusion we include several numerical examples in order to demonstrate the efficiency of proposed numerical methods.}

\end{abstract}


\end{frontmatter}

% \linenumbers

%% mai

\section*{Introduction}

%� ���������� ������� ����������� ���������� ������� ���������� ������� ������� ������������ ��������� ��������� I ����. ��� ������ ����� ��������� �� ��������� �������: �������� � ���������� II ����, ���������� ���������������� ����������, ������������� ��� ��������� ����������, ������ ������������� � �.�.
In this article we continue our studies of the novel class of linear Volterra (evolutionary)
integral equations (VIE) of the first kind with piecewise continuous kernels. %The problem of solution
The solution
of linear integral equations of the first kind is of course classical problem and has been addressed
by numerious authors. But only few authors studied these equations in case of jump
discontious kernels. In general,
VIE of the first kind can be solved by reduction to equations of the second kind,  regularization algorithms developed for Fredholm equations can be also applied as well as direct discretization methods.

%����� ��������, ��� ������� ������������ ��������� ������� ���� ����� ���� �������������, � ������ ��������� ����� ��������� ��������� � ������������. ��� ������� � ���, ��� ��������, ����������� ����� ������ ���������, ���������� ��������������� ������������ ������� �� �� ��� ��� ������������, � �� ��� ����� �����. ������� �������� �������� �� ���������. ������� ���������� ��������� �������� ������� � �������� ������ ������ � ������ �������� (�������� ������ ������ ������ ���� � ����� ������� �������).

From the other hand, it is known that solutions of integral equations of the first kind can be unstable and this is a well known  ill-posed problem. This is due to the fact that the Volterra operator maps the considered solution space into its narrow part only. Therefore, the inverse operator is not bounded. It is necessary to assess the proximity of the solutions and the proximity of the
%right parts
right-hand side
using the different metrics. In addition, the proximity of the
%right parts
right-hand side
should be in a stronger metric.
%����� ����, ��� �������� � \cite{Sidorov-book}, ������� ��������� ���� \eqref{e2015-03} ����� ��������� ������������ ���������� � ���� ��������������� ��� \(t\to 0\).
Moreover, as shown in \cite{Sidorov-book}, solutions of the VIE can contain arbitrary constants and can be unlimited as \(t\to 0\). Here readers may also refer to \cite{Sidorov-FR, SidorovMarkovaAiT, SidorovRM, SidorovDU}, where  the problems of existence, uniqueness and asymptotic behavior of solutions of equations of this type are  explored.

Evolutionary integral equations are in the core of many mathematical models in
physics, economics and ecology.
Excellent historical overview of the results concerning the
VIEs of the first kind is given by H. Brunner  in the paper ``1896 -- 1996:
One hundred years of Volterra integral equations of the first kind'' \cite{brunner} .
The theory of integral models of evolving systems was initiated  in the  works of  L. Kantorovich, R. Solow and V. Glushkov in the mid-20th Century. Here readers may refer to the papers of \cite{kantor} and \cite{solow}.  Such theory  employs  the VIEs of the first kind where bounds of the integration interval can be  functions of time. Here readers may refer e.g. to the monograph  \cite{hriton}.
These models take into account
the memory of a dynamical system when its past impacts its future evolution.
The memory is implemented in the existing technological and financial structure of physical capital (equipment).
The memory duration is determined by the age of the oldest capital unit (e.g. equipment) still employed.

The paper \cite{Boikov-Tynda11} is devoted to the construction of iterative numerical algorithm for the systems of nonlinear Volterra-type equations 
related to the Vintage Capital Models (VCMs) \cite{hriton}:
\[
  \begin{cases}
     x(t)=\int\limits_{y(t)}^{t}H(t,\tau,x(\tau))d\tau,\\
     \int\limits_{y(t)}^{t}K(t,\tau,x(\tau))d\tau
     =f(t),
  \end{cases}\;t\in[t_0,T),\;t_0< T \leqslant \infty,
\]
with unknown functions \(x(t)\) and \(y(t)\) satisfying the initial
conditions: \(y(t_0)=Y_0<t_0,\; x(\tau)\equiv  \varphi_0(\tau),\;\tau\in(-\infty,t_0].\)

Numerical methods which are optimal with respect to complexity order were constructed in paper \cite{Tynda-18} for VIEs with certain weakly singular kernels.

First results in studies of the Volterra equations with discontinuous kernels were
formulated by G.C. Evans~\cite{Evans} in the beginning of XX century.
Results in the spectral theory of integral operators with discontinuous kernels
were obtained by A.P. Khromov in his paper~\cite{Khromov}.  Some results concerning the general approximation theory for integral equations with discontinuous kernels are presented in paper~\cite{Anselone}.

%��� ��������� ������� ������������ ������������ ��������� ��������� I ���� �������� ��������� ��������. ���� �� ��� ������� � ���������� ������������ ���������������� ����������, ������������� ��� ������������ ��������� ���������� I ����. ������ ��� ���� ������ �������� � ������� �������������� ������ ��������� � ������ ��������, ������������ ������ ������������ ��������� ��������� � ���������� ������������ ���������� �������������� ��������� ����������. ������ ������ ������� �� ���������� ������ ������������� �������� ���������. ����� �������� ���������� �������������� ������������� ������� � ������������ �������� ������. ����� ����������� � ���� ������ �������� ������������� ���������������� ������� ������ ������� �������������, ��� ���� ��������� ������������� ������ ��� �������������, ��������� � ������� �������� ������. ��� ���� ��� ������������� ���������� �������� ���� ������������ ������� ���������� ������� �������� (������� ��������������� � ��������). ������� �������-������, ������� � ��. (������� 2 � ����) ��������� ������������ ���������. ��������� �������� ���������������� ������ ��������� ���������� ��������� � ����� \cite{Verlan-Sizikov}.
There are several approaches available for numerical solution of Volterra integral equations of the first kind. One of them is to apply classical regularizing algorithms developed for Fredholm integral equations of the first kind. However, the problem reduces to solving algebraic systems of equations with a full matrix, an important advantage of the Volterra equation is lost and there is a significant increase in arithmetic complexity of the algorithms. The second approach is based on a direct discretization of the initial equations. Here one may face an instability of the approximate solution %to
because of errors in the initial data. The regularization properties of the direct discretization methods are optimal in this sense, where the discretization step is the regularization parameter associated with the error of the source data. However, only low-order quadrature formulas (midpoint quadrature or trapezoidal formulas) are suitable for approximation of the integrals. The Newton-Cotes formulas, Gregory and others (the second order and higher orders) generate divergent algorithms. The detailed description of regularizing direct numerical algorithms is described in the monograph \cite{Kythe}.

%���������� ��������, ��� ���������� ���� ���������� � ��������� \eqref{e2015-01} � ����� \eqref{e2015-03} ������ �������������� ��-�� �������� ���� \eqref{e2015-02}. ����� ��� ���������� ������������� ���������� ���������� ����� ��� ������ ����� ��������� \(N\) ������ �������� �� ����� ������� \(\alpha_i(t)\) �, ��������������, �� ����� ���� ��������� � ������������ �������� ������.
It should be noted that %the application of
it is very difficult to apply
these algorithms to solve the equation \eqref{e2015-01} in the form of \eqref{e2015-03} %is very difficult
because of %breaks of the kernels
the kernel discontinuities \eqref{e2015-02} as described in Section 1.
%������  - ������ �������������, ���������� �� �������-���������� � �������-�������� ������������� ������� ������� (������� � ������� ������� ��������, ��������������).
The adaptive mesh should depend on the curves of the jump discontinuity for each number \(N\) of divisions of the considered interval and therefore this mesh can not be linked to the errors in the source data. It is needed to correctly approximate the integrals.

%���� ������������ ��� ������� � ���������� ������� ��������� \eqref{e2015-03}.
%������ ������ ������� �� ��������������� ����������� ���� ��������� �������� ����������� ������� � ����������� ���������� ����������� ���������������� ������������ ���������.

Below we continue our studies \cite{Sidorov-book,vestnik} and propose two approaches for the numerical solution for Volterra integral equations
of the first kind
with piecewise continuous kernels.
The first approach is a direct discretization based on piecewise constant and piecewise linear approximations of the exact solution (the first and the second order of accuracy, respectively).
The second approach is based on the preliminary determination of the two acceleration values of the unknown function and then we employ the special regularizing iterative procedure.

The paper is organized as follows. In Section 1, we describe the
problem give some statements concerning the existence and uniqueness of solutions of VIEs with discontinuous kernels.
Section 2 is dedicated to direct discretization numerical methods based on the piecewise constant and piecewise linear approximation of the exact solution.
In Section 3, we describe the regularization method for linear first kind VIEs of this class.
The modified Newton-Kantorovich iterative process for nonlinear VIEs is suggested in Section 4.
The numerical examples  are given in Section 5.

\section{Problem statement}

%�������� ���������� ������������ � ������ ������ �������� ������������ ��������� I ���� ���������� ����
The object of our interest is the following integral equation of the first kind
\begin{equation}\label{e2015-01}
  \int\limits_{0}^{t}K(t,s)x(s)\;ds=f(t), \; t\in[0,T],
\end{equation}
%��� ���� \(K(t,s)\) ������ �������� ������� �� ������ \(\alpha_i(t),\;i=1,2,\ldots,n-1,\) � ����� ���
where the kernel \(K(t,s)\) is discontinuous along continuous curves \(\alpha_i(t),\;i=1,2,\ldots,n-1,\) and is of the form
\begin{equation}\label{e2015-02}
  K(t,s)=\left\{
           \begin{array}{ll}
             K_1(t,s), &  \alpha_0(t)< s < \alpha_1(t); \\
             K_2(t,s), &  \alpha_1(t)< s < \alpha_2(t); \\
             \cdots  \\
             K_n(t,s), &  \alpha_{n-1}(t)< s < \alpha_n(t). \\
           \end{array}
         \right.
\end{equation}
%�����
Here
$
 \alpha_0(t)\equiv 0,\; \alpha_0(t)<\alpha_1(t)<\ldots<\alpha_n(t)\equiv t,\; f(0)=0.
$
%�����������, ��� ���� \(K_i(t,s)\) � ������ ����� \(f(t)\) � ��������� \eqref{e2015-01} �������� ������������ � ���������� �������� ���������. ������� \(\alpha_i(t)\in C^1[0,T]\) � �������� ������������. ����� ����
Let us assume, that the kernels \(K_i(t,s)\) and the
%right part
right-hand side
\(f(t)\) in the equation \eqref{e2015-01} are continuous and sufficiently smooth functions. The functions \(\alpha_i(t)\in C^1[0,T]\) are not decrescent. Moreover
\[
  \alpha_1'(0)\le\alpha_2'(0)\le\ldots\le\alpha_{n-1}'(0)<1.
\]
%��������� ��������� \eqref{e2015-01} � ����������� ����
Let us rewrite the equation \eqref{e2015-01}
\begin{equation}\label{e2015-03}
  \int\limits_{0}^{\alpha_1(t)}K_1(t,s)x(s)\;ds+\int\limits_{\alpha_{1}(t)}^{\alpha_2(t)}K_2(t,s)x(s)\;ds+
  \cdots+\int\limits_{\alpha_{n-1}(t)}^{t}K_n(t,s)x(s)\;ds=f(t), \; t\in[0,T].
\end{equation}
It is to be noted that conventional Glushkov integral model of evolving systems is the special case of this equation  where all the functions $K_i(t,s)$
 are zeros except of $K_n(t,s)$.

%************************************************************************************************************
%\section{������ �������������}
%\subsection{�������-���������� �������������}
\section{Direct discretization}
\subsection{Piecewise constant approximation}

%��� ���������� ���������� ������� ��������� \eqref{e2015-03} �� ������� \([0,T]\) (� �������� ������������� ������������� ������������ �������) ������ ����� ����� (������������� �����������)
Let us introduce the mesh nodes (not necessarily uniform) to construct the numeric solution of the equation \eqref{e2015-03} on the interval \([0,T]\)
%(in conditions of existence of a unique continuous solution)
(if the unique continuous solution exists)
\begin{equation}\label{e2015-04}
  0=t_0<t_1<t_2<\ldots<t_N=T, \;\: h=\max\limits_{i=\overline{1,N}}(t_i-t_{i-1})=O(N^{-1}).
\end{equation}
The approximate solution is determined as the following piecewise constant function
\begin{equation}\label{e2015-05}
   x_N(t)=\sum\limits_{i=1}^{N}x_i\delta_i(t),\; t\in(0,T], \;
   \delta_i(t)=\left\{
            \begin{array}{ll}
              1, & \hbox{ } t\in \Delta_i=(t_{i-1},t_i]; \\
              0, & \hbox{ } t\notin\Delta_i
            \end{array}
          \right.
\end{equation}
with the undefined coefficients \(x_i,\;i=\overline{1,N}\).
We differentiate the both parts of the equation \eqref{e2015-03} with respect to \(t\) to determine \(x_0=x(0)\)
\[
 f'(t)=\sum\limits_{i=1}^n \Biggl(\int\limits_{\alpha_{i-1}(t)}^{\alpha_i(t)}
       \frac{\partial K_i(t,s)}{\partial t}x(s)\;ds+\alpha'_i(t)K_i(t,\alpha_i(t))x(\alpha_i(t))-
\]
\[
       -\alpha'_{i-1}(t)K_i(t,\alpha_{i-1}(t))x(\alpha_{i-1}(t))\Biggr).
\]
From the last expression we obtain
\begin{equation}\label{e2015-06}
   x_0=\frac{f'(0)}{\sum\limits_{i=1}^n K_i(0,0)\left[\alpha'_i(0)-\alpha'_{i-1}(0)\right]}.
\end{equation}
%����� ��������������, ��� ����������  ���������  \eqref{e2015-03} ������, ��� ����������� �  \eqref{e2015-06} �� ���������� � ����.
Here it is assumed,
%equations \eqref{e2015-03} are such
that the denominator of \eqref{e2015-06} must be not zero.
%������ ����� ����������� \(f_k=f(t_k), \;k=1,\ldots,N\). ��� ����������� ������������ \(x_1\) ������� �������� ��������� � ����� \(t=t_1\):
We introduce the denotation \(f_k=f(t_k), \;k=1,\ldots,N\) and write the initial equation in the point \(t=t_1\) to define the coefficient \(x_1\):
\begin{equation}\label{e2015-07}
  \sum\limits_{i=1}^n \int\limits_{\alpha_{i-1}(t_1)}^{\alpha_i(t_1)}K_i(t_1,s)x(s)\;ds=f_1.
\end{equation}
%��� ��� �� ������ ���� ����� ���� �������� �������������� \(\alpha_i(t_1)-\alpha_{i-1}(t_1)\) � \eqref{e2015-07} �� ����������� \(h\), � ������������ ������� ��������� �������� \(x_1\), ��, �������� ������������ ������� ������� ���������������, �����
{Since at this stage the lengths of all integration intervals \(\alpha_i(t_1)-\alpha_{i-1}(t_1)\) in \eqref{e2015-07} don't exceed \(h\) then based on the midpoint quadrature rule we have
\begin{equation}\label{e2015-08}
   x_1=\frac{f_1}{\sum\limits_{i=1}^n (\alpha_i(t_1)-\alpha_{i-1}(t_1))K_i(t_1,\frac{\alpha_i(t_1)+\alpha_{i-1}(t_1)}{2})}.
\end{equation}
%����������� ������, ��� ��� ������� �������� \(x_2,x_3,\ldots,x_{k-1}\).
%��������� ��������� \eqref{e2015-01} � ����
Let us suppose now that we have already found the values \(x_2,x_3,\ldots,x_{k-1}\).
We rewrite the equation \eqref{e2015-01} as
\begin{equation}\label{e2015-09}
  \int\limits_{t_{k-1}}^{t}K(t,s)x(s)\;ds=f(t)-\int\limits_{0}^{t_{k-1}}K(t,s)x_N(s)\;ds
\end{equation}
%� ��������� ���������� ���������� ��������� � ����� \(t=t_k\):
and we require that the last equality hol\;ds for the point \(t=t_k\)
\begin{equation}\label{e2015-10}
  \int\limits_{t_{k-1}}^{t_k}K(t_k,s)x(s)\;ds=f_k-\int\limits_{0}^{t_{k-1}}K(t_k,s)x_N(s)\;ds.
\end{equation}
%� ������ \eqref{e2015-05} �����
Taking into account \eqref{e2015-05} we have
\[
   x_k\int\limits_{t_{k-1}}^{t_k}K(t_k,s)\;ds=f_k-\sum\limits_{j=1}^{k-1} x_j \int\limits_{t_{j-1}}^{t_{j}}K(t_k,s)\;ds.
\]
%������
Thus
\begin{equation}\label{e2015-11}
 x_k=\frac{f_k-\sum\limits_{j=1}^{k-1} x_j \int\limits_{t_{j-1}}^{t_{j}}K(t_k,s)\;ds}{\int\limits_{t_{k-1}}^{t_k}K(t_k,s)\;ds}.
\end{equation}
%��� ���� ��� ���������� ���������� ���� \(\int\limits_{t_{j-1}}^{t_{j}}K(t_k,s)\;ds\) � \eqref{e2015-11} ������������ ��������� ������� ������� ���������������, ����������� �� ��������������� ����� �����, ����������� ��� ������ ���������� �������� \(N\) � ������ \(\alpha_i(t)\) ��������  ���� \(K(t,s)\).
Herewith we calculate the integrals of the form \(\int\limits_{t_{j-1}}^{t_{j}}K(t_k,s)\;ds\) in \eqref{e2015-11} using the midpoint quadrature formulas with auxiliary mesh nodes related to the curves \(\alpha_i(t)\) of the kernels \(K(t,s)\) for each value of \(N\).
%�������� ��������, ��� ����������� ������ ����� ���:
It is easy to notice that the error of the method is
\begin{equation}\label{e2015-12}
 \varepsilon_N=\|x(t)-x_N(t)\|_{C_{[0,T]}}=O\left(\frac{1}{N}\right).
\end{equation}

%****************************************************
%****************************************************
%****************************************************
%\subsection{�������-�������� �������������}
\subsection{Piecewise linear approximation}
%����� ������ ������������ ������� ������������ �������-�������� ������� ����
We suppose that the approximate solution is a piecewise linear function of the following form
\begin{equation}\label{e2015-13}
   x_N(t)=\sum\limits_{i=1}^{N}\left(x_{i-1}+\frac{x_i-x_{i-1}}{t_i-t_{i-1}}(t-t_{i-1})\right)\delta_i(t),\; t\in(0,T], \;
\end{equation}
%���, ��� � �����,
where
\[
   \delta_i(t)=\left\{
            \begin{array}{ll}
              1, & \hbox{for } t\in \Delta_i=(t_{i-1},t_i]; \\
              0, & \hbox{for } t\notin\Delta_i.
            \end{array}
          \right.
\]
%������������ \(x_i,\;i=\overline{1,N},\) ������������� ������� �������� �����������.
%��������� �� ������� \eqref{e2015-06} ����������� \(x_0\)  � �������� ��������� \eqref{e2015-10}, �������
We need to determine the coefficients \(x_i,\;i=\overline{1,N},\) of the approximate solution.
Determining by the \eqref{e2015-06} the coefficient \(x_0\) and taking into account the equality \eqref{e2015-10} we obtain
\[
  \int\limits_{t_{k-1}}^{t_k}\left(x_{k-1}+\frac{x_k-x_{k-1}}{t_k-t_{k-1}}(s-t_{k-1})\right)K(t_k,s)\;ds=
\]
\[
  =f_k-\sum\limits_{j=1}^{k-1}\int\limits_{t_{j-1}}^{t_{j}}\left(x_{j-1}+\frac{x_j-x_{j-1}}{t_j-t_{j-1}}(s-t_{j-1})\right)K(t_k,s)\;ds.
\]
%������, �������� \(x_k\), �����
Thus excluding \(x_k\) we have
\begin{equation}\label{e2015-14}
 x_k=x_{k-1}+\frac{f_k-x_{k-1}\int\limits_{t_{k-1}}^{t_k}K(t_k,s)\;ds-
  \sum\limits_{j=1}^{k-1}\left(x_{j-1}\int\limits_{t_{j-1}}^{t_j}K(t_k,s)\;ds+\frac{x_j-x_{j-1}}{t_j-t_{j-1}}\int\limits_{t_{j-1}}^{t_j}(s-t_{j-1})K(t_k,s)\;ds \right)}
 {\frac{1}{t_k-t_{k-1}}\int\limits_{t_{k-1}}^{t_k}(s-t_{k-1})K(t_k,s)\;ds},
\end{equation}
%���
 where \( k=1,2,\ldots,N.\)

%��� ���������, �������� � \eqref{e2015-14}, ���������������� � ������� ������������ ������ ������� ���������������, ����������� �� ��������������� ����� ����� (��� ����� �������� ������� \(\alpha_i(t_j)\) �������� ������������� ��������� ����� ���� ����� ��� ������ ���������� �������� \(N\)).
We approximate the integrals in \eqref{e2015-14} by using the midpoint quadrature formulas based on auxiliary mesh nodes so that the values of the functions \(\alpha_i(t_j)\) are a subset of the set of this mesh points at each particular value of \(N\).
%����������� ������ ��� ����� ������������� ����� ���:
The error of this approximation method is
\begin{equation}\label{e2015-15}
 \varepsilon_N=\|x(t)-x_N(t)\|_{C_{[0,T]}}=O\left(\frac{1}{N^2}\right).
\end{equation}

%************************************************************************************************************
%\section{������������ �����}
\section{Iterative method}

%����� ������������ ���� \(K_i(t,s)\) � ��������������� �������� ����������� �������� ��������������� ���������, �.�.
Let the kernels \(K_i(t,s)\) be a symmetric functions in their domains, i.e.
\[
  K_i(t,s)=K_i(s,t), \; i=1,2,\ldots, n.
\]
%������������ ������� ��������� \eqref{e2015-01} �� ����� \eqref{e2015-04} ����� ������ � ���� �������-���������� ������� \eqref{e2015-05}. ��� �����
%���������, ��� � �����, ��������� �������� \(x_0\) � \(x_1\), ��������� ������� \eqref{e2015-06}, \eqref{e2015-08}.
%��������� ������ ��������� \eqref{e2015-01} � ����
We search the approximation solution of the equation \eqref{e2015-01} at the mesh \eqref{e2015-04} as a piecewise constant function like \eqref{e2015-05}. To do this
we define initial values of the \(x_0\) and \(x_1\) with the formulas \eqref{e2015-06}, \eqref{e2015-08}.
We rewrite the equation \eqref{e2015-01}
\begin{equation}\label{e2015-21}
  \int\limits_{t_1}^{t}K(t,s)x(s)\;ds=f(t)-\int\limits_{0}^{t_1}K(t,s)x_N(s)\;ds
\end{equation}
%� ��������� \(g(t)=f(t)-\int\limits_{0}^{t_1}K(t,s)x_N(s)\;ds\).
and designate \(g(t)=f(t)-\int\limits_{0}^{t_1}K(t,s)x_N(s)\;ds\).
%��� ����������� �������� \(x_k\) �������� ������������� ������� \eqref{e2015-05} ���������� ��������� ������������ �������:
To define the values \(x_k\) of the required approximation solution \eqref{e2015-05} we use the following iterative process:
\begin{equation}\label{e2015-22}
  x^{(m+1)}(t)=x^{(m)}(t)+\gamma\left(g(t)-\int\limits_{t_1}^{t}K(t,s)x^{(m)}(s)\;ds\right),
  \; m=0,1,\ldots,
\end{equation}
%��� \(\gamma\) --- ������������� �������� �������������, \(m\) --- ����� ��������.
where \(\gamma\) is a positive regularization parameter and \(m\) is a number of the iteration.

%��������� ����������� \(x^{(0)}(t)\) ������������ �� ��������� ������ � ������ ������� (��� �� �������) ��� ���������� \(x^{(0)}(t)\equiv g(t)\).
We define the initial approximate value \(x^{(0)}(t)\) from the aprioristic data (if we have it) of the exact solution or we suppose \(x^{(0)}(t)\equiv g(t)\).
%��������, ��� ���� ���������� ����� ������� �������������� ������������������ \(x^{(m)}(t)\) �������� � ��������� ������� \(\tilde{x}_{\gamma}(t)\), �� ��� ������� ������������� ��������� \eqref{e2015-22} ��� ����� \(\gamma\ne 0\).
Obviously, if the functional sequence \(x^{(m)}(t)\) converge to a function \(\tilde{x}_{\gamma}(t)\) then that function satisfy \eqref{e2015-22} for all \(\gamma\ne 0\).
%�������� \(x_k, \; k=2,3,\ldots,N,\) ������������ ����� ��������������� �� �����
The values \(x_k, \; k=2,3,\ldots,N,\) can be defined successively as
\begin{equation}\label{e2015-23}
  x^{(m+1)}_k=x^{(m)}_k+\gamma\left(g(t_k)-\int\limits_{t_1}^{t_k}K(t_k,s)x^{(m)}_N(s)\;ds\right),
  \; m=0,1,\ldots.
\end{equation}
%��� ���� ��� ���������� ���������� � \eqref{e2015-23} ������������ ��������� ������� ������� ��������������� ��� ��������, ����������� �� ��������������� ����� �����, ����������� ��� ������ ���������� �������� \(N\) � ������ \(\alpha_i(t)\) ��������  ���� \(K(t,s)\).
%����������� �������� ��������� ������������� \(\gamma\) ���������� ������������ �������� �� ������� ������������� �������
Herewith to calculate the integrals in the \eqref{e2015-23} we use midpoint quadrature or trapezoidal formulas based on auxiliary mesh nodes related to the curves \(\alpha_i(t)\) of the kernels \(K(t,s)\) for each value of \(N\).
In practice, we choose the optimal value of the regularization parameter \(\gamma\) with the following condition
\begin{equation}\label{e2015-24}
  \varepsilon_N^{(m)}=\left\|\sum\limits_{i=1}^n \int\limits_{\alpha_{i-1}(t)}^{\alpha_i(t)}K_i(t,s)x^{(m)}_N(s)\;ds-f(t)\right\|_{C_{[0,T]}}\to \min
\end{equation}
%��� ���������� ������� \(m\).
for large enough \(m\).

%-------------------------------------------------------------------------
%      Nonlinear equations
%-------------------------------------------------------------------------
\section{Nonlinear equations}
In this section we consider the extension of \eqref{e2015-01} to the case of nonlinear dependency \(K(t,s,x(s))\):
\begin{equation}\label{e2015-31}
  \int\limits_0^t K(t,s,x(s))\;ds = f(t), \quad t\in[0,T], \;f(0)=0,
\end{equation}
where
\begin{equation}\label{e2015-32}
    K(t,s,x(s)) = \left\{ \begin{array}{ll}
         \mbox{$K_1(t,s)G_1(s,x(s)), \,\, t,s \in m_1$}, \\
         \mbox{\,\, \dots \,\, \dots \dots \dots } \\
         \mbox{$K_n(t,s)G_n(s,x(s)), \,\, t,s \in m_n$}. \\
        \end{array} \right.
\end{equation}
Here \( m_i = \{ t,s  \bigl | \alpha_{i-1}(t) < s \leqslant \alpha_i(t) \},\)   \( \alpha_0(t)=0,\; \alpha_n(t)=t,\; i=\overline{1,n}\),
The functions \(K_i\), \(f(t)\),  \(\alpha_i(t)\) have continuous derivatives with respect to  $t$ at $t,s \in \overline{m_i},$  $ K_n(t,t) \neq 0,$   $\alpha_i(0)=0,$ $\,\,\,\,0 < \alpha_1(t)<\alpha_2(t)< \cdots < \alpha_{n-1}(t)<t$.
The functions $\alpha_1(t), \dots , \alpha_{n-1}(t) $  should increase in small neighbourhood $0\leq t \leq \tau$ at least.

The following theorem states the existence and uniqueness conditions of solution of equation \eqref{e2015-31}. The proof is similar with proof of the Theorem 3.2 in 
the monograph \cite{Sidorov-book}.

\begin{theorem} \label{th1nlin} %{the existence and uniqueness conditions of solution of equation}
Let for $t\in[0,T]$ the following conditions takes place:
 $K_i(t,s), G_i(s,x(s))$ are continuous,
 $\, i=\overline{1,n}$,  $\alpha_{i}(t)$ and $f(t)$
have continuous derivatives for $t$,  $K_n(t,t)\neq 0$,
$0=\alpha_0(t)<\alpha_1(t)< \dots < \alpha_{n-1}(t)<\alpha_n(t)=t$
 for $t \in (0,T]$, $\alpha_i(0)=0$, $f(0)=0.$
Let the functions $G_i(s,x(s))$ satisfy Lipschitz condition
$|G_i(s,x_1(s))-G_i(s,x_2(s))| \leq q_i|x_1 - x_2|$, $\forall x_1, x_2 \in \mathbb{R}^1$,
$q_n + \sum\limits_{i=1}{n-1}{\alpha_i^{\prime}(0)|K_n(0,0)^{-1} (K_i(0,0)-K_{i+1}(0,0))| (1+q_i)<1}.$
Then $\exists \tau>0$ such that the equation (\ref{e2015-31}) has an unique local solution in $\mathcal{C}_{[0,\tau]}.$
Furthermore, if {$\min\limits_{\tau \leq t \leq T} (t-\alpha_{n-1}(t)) = h>0$} then the the solution can be constructed on $[\tau, T]$
using the step method combined with successive approximations. Thereby the equation (\ref{e2015-31}) has the unique global solution in ${\mathcal{C}_{[0,T]}}$.
\end{theorem}

\subsection{Linearization}

In order to approximate solution of \eqref{e2015-31} we introduce the nonlinear integral operator
\begin{equation}\label{e2015-33}
  (Fx)(t)\equiv \sum\limits_{i=1}^{n}\int\limits_{\alpha_{i-1}(t)}^{\alpha_i(t)} K_i(t,s)G_i(s,x(s))\;ds - f(t).
\end{equation}
The equation \eqref{e2015-31} can be written in an operator form as follows:
\begin{equation}\label{e2015-34}
  (Fx)(t)=0.
\end{equation}
In order to construct an iterative numerical method to equation \eqref{e2015-31}, we first linearize the operator \eqref{e2015-33} according to a
modified Newton-Kantorovich scheme \cite{Kantorovich2}:
\begin{equation}\label{e2015-35}
     x_{m+1}=x_m-[F'(x_0)]^{-1}(F(x_m)) ,\; m=0,1,\ldots,
\end{equation}
where \(x_0(t)\) is the initial approximation. Then, the approximate solution of \eqref{e2015-34} could be determined as the following limit of sequence:
\begin{equation}\label{e2015-36}
  x(t)=\lim\limits_{m\to\infty}x_m(t).
\end{equation}
The derivative \(F'(x_0)\) of the nonlinear operator \(F\) at the point \(x_0\) is defined as follows:
\[
  F'(x_0)=\lim_{\omega\to 0}\frac{F(x_0+\omega x)-F(x_0)}{\omega}=
\]
\[
  =\lim_{\omega\to 0}\frac{1}{\omega}\sum\limits_{i=1}^{n}\int\limits_{\alpha_{i-1}(t)}^{\alpha_i(t)}
  K_i(t,s)\left[G_i(s,x_0(s)+\omega x(s))-G_i(s,x_0(s))\right]\;ds.
\]
Implementing the limit transition under the integral sign, we finally get:
\begin{equation}\label{e2015-37}
    F'(x_0)(t)=\sum\limits_{i=1}^{n}\int\limits_{\alpha_{i-1}(t)}^{\alpha_i(t)}
    K_i(t,s)G_{ix}(s,x_0(s))x(s)\;ds, \text{ where } G_{ix}(s,x_0(s))=\left.\frac{\partial G_i(s,x(s))}{\partial x}\right|_{x=x_0}.
\end{equation}
%%-----------------------
Thus, we obtain the operator form of Newton-Kantorovich scheme as follows:
\begin{equation}\label{e2015-38}
    F'(x_0(t))\Delta x_{m+1}(t)=-F(x_m), \; \Delta x_{m+1}=x_{m+1}-x_m,
\end{equation}
or in the extended form
\[
\sum\limits_{i=1}^{n}\int\limits_{\alpha_{i-1}(t)}^{\alpha_i(t)}
    K_i(t,s)G_{ix}(s,x_0(s))\Delta x_{m+1}(s)\;ds=
    f(t)-\sum\limits_{i=1}^{n}\int\limits_{\alpha_{i-1}(t)}^{\alpha_i(t)}
    K_i(t,s)G_{i}(s,x_m(s))\;ds.
\]
%%-----------------------
We rewrite the last equation as follows
\begin{equation}\label{e2015-39}
    \sum\limits_{i=1}^{n}\int\limits_{\alpha_{i-1}(t)}^{\alpha_i(t)}
    K_i(t,s)G_{ix}(s,x_0(s)) x_{m+1}(s)\;ds=\Psi_m(t), \; m=0,1,2,\ldots,
\end{equation}
where
\[
  \Psi_m(t)=f(t)+\sum\limits_{i=1}^{n}\int\limits_{\alpha_{i-1}(t)}^{\alpha_i(t)}
   K_i(t,s)\left[G_{ix}(s,x_0(s)) x_{m}(s)-G_{i}(s,x_m(s))\right] \;ds.
\]
Equations \eqref{e2015-39} are now linear Volterra equations of the first kind with respect to the unknown function \(x_{m+1}(t)\).
Note that the kernels \(K_i(t,s)G_{ix}(s,x_0(s))\), \(i=\overline{1,n},\)   remain constant during each iteration \(m\).
Since equations \eqref{e2015-39} have the form \eqref{e2015-03} we can apply the methods suggested in Section 2 and Section 3 to solve them numerically.
Thus, solving the equations \eqref{e2015-39}, we get a sequence of approximate functions \(x_{m+1}(t)\). And then, using formula \eqref{e2015-36},
we obtain the approximate solution of \eqref{e2015-31} with an accuracy depending on \(m\).

\subsection{The convergence theorem }
Let \(C[0,T]\) be a Banach space of continuous functions equipped
with the standard norm \( \|x\|_{C[0,T]}=\max\limits_{t\in [0,T]}|x(t)| \).
The following theorem of convergence (based on the general theory proposed in the classical monograph
\cite{Kantorovich2}) for iterative process \eqref{e2015-39} takes
place:
%%------------------
\begin{theorem}
   Let the operator \(F\) has a continuous second derivative in the sphere
   \(\Omega_0\;(\|x-x_0\|\leqslant r)\) and the following conditions
   hold:
   \begin{enumerate}
      \item Equation \eqref{e2015-39} has a unique solution in
      \([0,T]\) for \(m=0\), i.e. there exists \(\Upsilon_0=[F'(x_0)]^{-1}\);
      \item \(\|\Delta x_1\|\leqslant\eta; \)
      \item \(\|\Upsilon_0F''(x)\|\leqslant L,\;\;x\in\Omega_0\).
   \end{enumerate}
If also
\(
  h=L\eta<\frac12 \text{ and }
  \frac{1-\sqrt{1-2h}}{h}\eta\leqslant
  r\leqslant\frac{1+\sqrt{1-2h}}{h}\eta,
\)
then equation \eqref{e2015-31} has a unique solution \(x^*\) in
\(\Omega_0\), process \eqref{e2015-39} converges to \(x^*\), and the
velocity of convergence is estimated by the inequality
  \[\|x^*-x_m\|\leqslant \frac{\eta}{h}(1-\sqrt{1-2h})^{m+1},
  \;m=0,1,\ldots.\]
\end{theorem}
%%------------------
In order to prove this theorem we show that equation
\eqref{e2015-39} is uniquely solvable (including the case \(m=0\)),
i.e. condition 1 of the theorem holds. Then we verify the
boundedness of the second derivative \(\bigl[F''(x_0)\bigr](x)\))
for estimating the constant \(L\) in condition 3.

One can verify that the necessary condition for the second
derivative \(\bigl[F''(x_0)\bigr](x)\) to be bounded is a
differentiability of the initial approximation \(x_0(t)\) as well as
the functions \(K_i\) with respect to second variable.

%%------------------------------------------------------------------------------------
%%------------------------------------------------------------------------------------
\subsection{Generalized numerical method for nonlinear equations}
In this section we offer for nonlinear weakly regular Volterra equations common numerical method based on using midpoint quadrature rule.

To find numerical solution of the equation (\ref{e2015-31}) on the interval \([0,T]\) we introduce the following mesh  (the mesh can be non-uniform)
\begin{equation}\label{e2015-40}
  0=t_0<t_1<t_2<\ldots<t_N=T, \;\: h=\max\limits_{i=\overline{1,N}}(t_i-t_{i-1})=\mathcal{O}(N^{-1}).
\end{equation}
Let us search for the approximate solution of the equation \ref{e2015-31} as following piecewise constant function
\begin{equation}\label{e2015-41}
   x_N(t)=\sum_{i=1}^{N}x_i\delta_i(t),\; t\in(0,T], \;
   \delta_i(t)=\left\{
            \begin{array}{ll}
              1, & \hbox{for } t\in \Delta_i=(t_{i-1},t_i]; \\
              0, & \hbox{for } t\notin\Delta_i
            \end{array}
          \right.
\end{equation}
with  coefficients  \(x_i,\;i=\overline{1,N}\) are under determination.
In order to find \(x_0=x(0)\) we differentiate both sides of the equation \ref{e2015-31} with respect to \(t\):
\[
f'(t)=\sum\limits_{i=1}^n \Biggl(\,\,\,\int\limits_{\alpha_{i-1}(t)}^{\alpha_i(t)}
       \frac{\partial K_i(t,s)}{\partial t}G_i(s, x(s))  \; ds+\alpha'_i(t)K_i(t,\alpha_i(t))G_i(\alpha_i(t), x(\alpha_i(t)))-\]
     \[ - \alpha'_{i-1}(t)K_i(t,\alpha_{i-1}(t))G_i(\alpha_{i-1}(t),x(\alpha_{i-1}(t)))\Biggr).
\]
Thereby
\[
 f'(0)=\sum\limits_{i=1}^n \Biggl(\,\,\,\int\limits_{0}^{0}
      \frac{\partial K_i(0,0)}{\partial t}G_i(0, x(0))  \; ds+\alpha'_i(0)K_i(0,0)G_i(0, x(0))- \alpha'_{i-1}(0)K_i(0,0)G_i(0,x(0))\Biggr).
\]
In the last expression the coefficient \(x_0\) appears in the case of nonlinear dependency.
To find the coefficient \(x_0\) we use Van Wijngaarden\(-\)Dekker\(-\)Brent method. The implementation of this method are considered in detail in \cite{Recipes}.
Let us make the notation \(f_k=f(t_k), \;k=1,\ldots,N\).
The mesh point of the mesh \ref{e2015-40} which coincide with
  \(\alpha_i(t_j)\)  we still denote as \(v_{ij}\), i.e.  \(\alpha_{i}(t_j)\in\Delta_{v_{ij}}\). Obviously \(v_{ij}<j\) for \(i=\overline{0,n-1}\), \(j=\overline{1,N}\).
{It is to be noted that \(\alpha_i(t_j)\) are not always
coincide  with any mesh point. Here \(v_{ij}\)
is used as index of the segment
\(\Delta_{v_{ij}}\),  such as  \(\alpha_i(t_j) \in \Delta_{v_{ij}}\) (or its right-hand side).}
Let us now assume the coefficients \(x_0,x_1,\ldots,x_{k-1}\) be known.
Equation \ref{e2015-31} defined in \(t=t_k\) as
\[
  \sum\limits_{i=1}^n \int\limits_{\alpha_{i-1}(t_k)}^{\alpha_i(t_k)}K_i(t_k,s)G_i(s,x(s)) \; ds=f_k,
\]
we can rewrite as follows:
$I_1(t_k)+I_2(t_k)+\cdots + I_n(t_k) = f_k, $
where
$$I_1(t_k)= \sum\limits_{j=1}^{v_{1,k}-1} \int\limits_{t_{j-1}}^{t_j} K_1(t_k,s) G_1(s,x(s))\,ds+
\int\limits_{t_{v_{1,k}}-1}^{\alpha_1(t_k)} K_1(t_k,s) G_1(s,x(s))\,ds,$$
$$I_n(t_k)= \int\limits_{\alpha_{n-1}(t_k)}^{t_{v_{n-1,k}}} K_n(t_k,s) G_n(s,x(s))\,ds +
\sum\limits_{j=v_{n-1,k}+1}^k \int\limits_{t_{j-1}}^{t_j} K_n(t_k,s) G_n(s,x(s))\,ds.$$
\begin{enumerate}
\item If  $v_{p-1,k} \neq v_{p,k},\, p=2,\dots, n-1$, then
$$I_p(t_k) = \int\limits_{\alpha_{p-1}(t_k)}^{t_{v_{p-1,k}}}K_p(t_k,s)G_p(s,x(s))\, ds +
\sum\limits_{j=v_{p-1,k}+1}^{v_{p,k}-1} \int\limits_{t_{j-1}}^{t_j} K_p(t_k,s) G_p(s,x(s))\,ds+\int\limits_{t_{v_{p,k}}-1}^{\alpha_p(t_k)} K_p(t_k,s) G_p(s,x(s))\, ds. $$
\item If $v_{p-1,k}=v_{p,k},\, p=2,\dots, n-1$, then
$$I_p(t_k)= \int\limits_{\alpha_{p-1}(t_k)}^{\alpha_p(t_k)} K_p(t_k,s) G_p(s,x(s))\,ds.$$
\end{enumerate}
The number of terms in each line of the last formula depends on an array \(v_{ij}\),
defined using the input data: functions \(\alpha_i(t),\;i=\overline{1,n-1}\),
and fixed (for specific $N$) mesh.
Each integral term we approximate using the midpoint quadrature rule, e.g.
\[
  \int\limits_{t_{v_{p,k}-1}}^{\alpha_p(t_k)}K_p(t_k,s)G_p(s,x(s)) \; ds\approx \left(\alpha_p(t_k)-t_{v_{p,k}-1}\right)K_p\left(t_k,\frac{\alpha_p(t_k)+t_{v_{p,k}-1}}{2}\right)G_p\left(\frac{\alpha_p(t_k)+t_{v_{p,k}-1}}{2}, x_N\left(\frac{\alpha_p(t_k)+t_{v_{p,k}-1}}{2}\right)\right).
\]
Moreover,
on those intervals where the desired function has been already determined,
we select \(x_N(t)\) (i.e. \(t\leqslant t_{k-1}\)).
On the rest of the intervales an unknown value $x_k$ appears in the last terms.
We explicitly define it and proceed in the loop for \(k \). The number of these terms is determined from the initial data \(v_{ij}\) analysis.
To find the coefficient \(x_0\) we also use Van Wijngaarden\(-\)Dekker\(-\)Brent method.
The maximum pointwise error of proposed numerical method
$\varepsilon^N=\max\limits_{0\leq i \leq N}|\bar{x}(t_i)-x^h(t_i)| $
%\end{equation}
has order of \({\mathcal O}\left(\frac{1}{N}\right)\).

%----------------------------

\section{Numerical examples}
\subsection{Linear equations}

Let us consider the following three problems using the uniform meshes only.
We define on \([0,T]\) the two-mesh difference $D^N$ for the $h=\frac{1}{N}$ and the $h=\frac{1}{2N}$
and the order of convergence $p^N$ based on the $D^N$ as follows
\begin{equation}\label{DN_eq}
D_N=\max\limits_{0\leq i \leq N, 0\leq j \leq 2N}|x^h_N(t_i) - x^h_{2N}(t_{j})|,
\end{equation}
with $t_i = t_{j}, i=2j,$
\begin{equation}\label{RN_eq}
p_N=\log_2{\frac{D_N}{D_{2N}}}.
\end{equation}
We use the $D^N$ and the $p^N$ to estimate the order of convergence for plorlems with unknown exact solutions.

Let us first address the equation
%\begin{equation}
$$\int\limits_0^{t/3} (1+t+s) x(s)\;ds-\int\limits_{t/3}^{t} x(s)\;ds =\frac{(2t+1)^\frac{3} {2}}{3} + \frac{\sqrt{3}(2t+3)^\frac{5} {2}}{45} - \frac{7 t^2}{18} - \frac{4}{15},\; t\in [0,\,2],
$$%\end{equation}
\begin{figure}[htb!] %T2
%{Computed maximum pointwise errors $\varepsilon^N,$  two-mesh difference $D^N$ and order of convergence $p^N$ in double precision arithmetic applied to problem 1 for various values of $h=1/N$}
\begin{flushright}
%{Tab. 4.1 Estimations for the 1st example}\vspace{-2mm}
{Tab. 5.1.1 The errors for various stepsizes $h$ of the 1st example}\vspace{-2mm}
\end{flushright}
\begin{center}
\begin{tabular}{|c|c|c|c|c|c|c|c|c|}	
\hline
\multicolumn{9}{|c|}{Piecewise constant approximation} \\ \hline
		$h$			& $1/32$   & $1/64$   & $1/128$  & $1/256$  & $1/512$  & $1/1024$ & $1/2048$ & $1/4096$ \\ \hline
\textbf{$\varepsilon$} 	& 0.097245 & 0.037330 & 0.020360 & 0.013031 & 0.005846 & 0.003016 & 0.001540 & 0.000781 \\ \hline
\textbf{${D_N}$} 		& 0.078452 & 0.027076 & 0.018686 & 0.009008 & 0.003580 & 0.002179 & 0.000852 & 0.000441 \\ \hline
\textbf{${p_N}$}		& 1.534798 & 0.535009 & 1.052673 & 1.330913 & 0.716153 & 1,354821 & 0,949923 &  - 		\\ \hline
\hline
\multicolumn{9}{|c|}{Piecewise linear approximation} \\ \hline
		$h$			& $1/32$   	& $1/64$   	& $1/128$  	& $1/256$  	& $1/512$  	& $1/1024$ 	& $1/2048$ 	& $1/4096$ 	\\ \hline
\textbf{$\varepsilon$} 	& 9.9445E-4 & 3.7228E-4 & 1.1005E-4 & 2.4769E-5 & 7.6136E-6 & 1.7694E-6 & 4.8759E-7 & 1.4031E-7 \\ \hline
\textbf{${D_N}$} 		& 9.3805E-4 & 3.3986E-4 & 1.0325E-4 & 2.1979E-5 & 6.7433E-6 & 1.6615E-6 & 4.3915E-7 & 1.2764E-7 \\ \hline
\textbf{${p_N}$}		& 1.464704 	& 1.718816 	& 2.231901 	& 1.704639 	& 2.020961 	& 1.919694 	& 1.782557 	&  - 		\\ \hline
\hline
\multicolumn{9}{|c|}{Iterative method} \\ \hline
		$h$			& $1/32$   & $1/64$   & $1/128$  & $1/256$  & $1/512$  & $1/1024$ & $1/2048$ & $1/4096$ \\ \hline
\textbf{$\varepsilon$} 	& 0.114115 & 0.044888 & 0.022213 & 0.013809 & 0.006252 & 0.003126 & 0.001475 & 0.000766 \\ \hline
\textbf{${D_N}$} 		& 0.088355 & 0.031906 & 0.019514 & 0.009391 & 0.003886 & 0,002224 & 0,000898 & 0,000448 \\ \hline
\textbf{${p_N}$}		& 1.469484 & 0.709318 & 1.055158 & 1.272992 & 0.805129 & 1.308369 & 1.003216 &  - 		\\ \hline
\end{tabular}
\end{center}
\label{ms_tab1}
\end{figure}
with known solution $\bar{x}(t)=\sqrt{2t+1}-1$.
Tab. \ref{ms_tab1}.1 shows computed maximum pointwise errors $\varepsilon^N=\max\limits_{0\leq i \leq N}|\bar{x}(t_i)-x^h(t_i)|,$ the two-mesh difference $D^N$
and the order of convergence $p^N$ in double precision arithmetic applied to problem 1 for various values of $h$.

Let us now consider the equation
%\begin{equation}
$$\int\limits_0^{t/8} (1-t \cdot s) x(s)\;ds + \int\limits_{t/8}^{3t/8} (t+s) x(s)\;ds - \int\limits_{3t/8}^{t} x(s)\;ds
=-\frac{t^5}{16384} + \frac{67t^4}{3072} - \frac{121t^3}{384},\; t\in [0,\,2],
$$%\end{equation}
where $\bar{x}(t)=t^2$ is exact solution.
Tab. \ref{ms_tab2}.2 shows $\varepsilon^N$, $D^N$, $p^N$ for various values of $h$.
\begin{figure}[htb!] %T2
%\tblcaptionnotes{Computed maximum pointwise errors $\varepsilon^N,$  two-mesh difference $D^N$ and order of convergence $p^N$ in double precision arithmetic applied to problem 2 for various values of N}
\begin{flushright}
%{Tab. 4.2 Estimations for the 2nd example}\vspace{-2mm}
{Tab. 5.1.2 The errors for various stepsizes $h$ of the 2nd example}\vspace{-2mm}
\end{flushright}
\begin{center}
\begin{tabular}{|c|c|c|c|c|c|c|c|c|}	
\hline
\multicolumn{9}{|c|}{Piecewise constant approximation} \\ \hline
				$h$			& $1/32$   & $1/64$   & $1/128$  & $1/256$  & $1/512$  & $1/1024$ & $1/2048$ & $1/4096$ \\ \hline
\textbf{$\varepsilon$} 	& 0.152263 & 0.084389 & 0.043314 & 0.021544 & 0.011043 & 0.005522 & 0.002759 & 0.001385 \\ \hline
\textbf{${D_N}$} 		& 0.073562 & 0.041074 & 0.021769 & 0.010589 & 0.005520 & 0.002805 & 0.001423 & 0.000715 \\ \hline
\textbf{${p_N}$}		& 0.840735 & 0.915910 & 1.039741 & 0.939704 & 0.976789 & 0.978239 & 0.992681 &  - 		\\ \hline
\hline
\multicolumn{9}{|c|}{Piecewise linear approximation} \\ \hline
				$h$			& $1/32$	& $1/64$   	& $1/128$  	& $1/256$  	& $1/512$  	& $1/1024$ 	& $1/2048$ 	& $1/4096$ 	\\ \hline
\textbf{$\varepsilon$} 	& 2.4351E-3 & 8.3747E-4 & 1.7396E-4 & 8.0309E-5 & 2.7166E-5 & 1.3062E-5 & 6.3391E-6 & 3.1053E-6 \\ \hline
\textbf{${D_N}$} 		& 2.4588E-3 & 8.3227E-4 & 1.6570E-4 & 6.4038E-5 & 2.9543E-5 & 1.5713E-5 & 8.2413E-6 & 4.3645E-6 \\ \hline
\textbf{${p_N}$}		& 1.562823 	& 2.328444 	& 1.371611 	& 1.116090 	& 0.910813 	& 0.931088 	& 0.917034 	&  - 		\\ \hline
\hline
\multicolumn{9}{|c|}{Iterative method} \\ \hline
			$h$				& $1/32$   & $1/64$   & $1/128$  & $1/256$  & $1/512$  & $1/1024$ & $1/2048$ & $1/4096$ \\ \hline
\textbf{$\varepsilon$} 	& 0.152264 & 0.084389 & 0.043315 & 0.021545 & 0.011043 & 0.005523 & 0.002760 & 0.001385 \\ \hline
\textbf{${D_N}$} 		& 0.073563 & 0.041074 & 0.021770 & 0.010589 & 0.005521 & 0.002805 & 0.001424 & 0.000716 \\ \hline
\textbf{${p_N}$}		& 0.840754 & 0.915884 & 1.039775 & 0.939564 & 0.976928 & 0.978051 & 0.991917 &  - 		\\ \hline
\end{tabular}
\end{center}
\label{ms_tab2}
\end{figure}

Finally we demonstrate the results obtained for the equation
%$$
%\begin{array}{c}
%%\begin{equation}
$$\int\limits_0^{t/8} (1+t+s) x(s)\;ds + \int\limits_{t/8}^{t/2} (2+t  s) x(s)\;ds + \int\limits_{t/8}^{3t/4} (t+s-1) x(s)\;ds -4\int\limits_{3t/4}^{t} x(s)\;ds=$$% \\
$$=\frac{1}{128}  \left(-4 - \frac{1}{8}  (16t + 69t^2 + 15t^3) - e^{\frac{t}{4}}  (t^2 - 13t + 12) + e^{t}  (4t^2 - 16t + 28) + e^{\frac{3t}{2}}  (14t + 20) - 32 e^{2t}\right),$$% \\
$$\; t\in [0,\,2],$$%
%%\end{equation}
%\end{array}
%\eqno(4.5)
%$$
with known solution $\bar{x}(t)=\frac{e^{2t}-1}{8}$.
Tab. \ref{ms_tab3}.3 shows $\varepsilon^N$, $D^N$, $p^N$ for various  $h$.

\begin{figure}[htb!] %T2
%\tblcaptionnotes{Computed maximum pointwise errors $\varepsilon^N,$  two-mesh difference $D^N$ and order of convergence $p^N$ in double precision arithmetic applied to problem 3 for various values of N}
%{\mbox{\tabcolsep=10pt\begin{tabular}{@{}cccc@{}}
\begin{flushright}
%{Tab. 4.3 Estimations for the 3rd example}\vspace{-2mm}
{Tab. 5.1.3 The errors for various stepsizes $h$ of the 3rd example}\vspace{-2mm}
\end{flushright}
\begin{center}
\begin{tabular}{|c|c|c|c|c|c|c|c|c|}	
\hline
\multicolumn{9}{|c|}{Piecewise constant approximation} \\ \hline
				$h$			& $1/32$   & $1/64$   & $1/128$  & $1/256$  & $1/512$  & $1/1024$ & $1/2048$ & $1/4096$ \\ \hline
\textbf{$\varepsilon$} 	& 0.449935 & 0.248336 & 0.127251 & 0.064714 & 0.033201 & 0.017015 & 0.008610 & 0.004350 \\ \hline
\textbf{${D_N}$} 		& 0.201599 & 0.121084 & 0.062537 & 0.031512 & 0.016186 & 0.008404 & 0.004260 & 0.002150 \\ \hline
\textbf{${p_N}$}		& 0.735482 & 0.953223 & 0.988773 & 0.961195 & 0.945483 & 0.980078 & 0.986308 &  - 		\\ \hline
\hline
\multicolumn{9}{|c|}{Piecewise linear approximation} \\ \hline
			$h$				& $1/32$	& $1/64$   	& $1/128$  	& $1/256$  	& $1/512$  	& $1/1024$ 	& $1/2048$ 	& $1/4096$ 	\\ \hline
\textbf{$\varepsilon$} 	& 2.1312E-2 & 8.5079E-3 & 4.6961E-3 & 2.0226E-3 & 1.8494E-3 & 7.2204E-4 & 3.8049E-4 & 1.5079E-4 \\ \hline
\textbf{${D_N}$} 		& 1.2804E-2	& 5.6290E-3 & 3.1706E-3 & 1.7218E-3 & 1.3342E-3 & 5.0412E-4 & 2.7926E-4 & 1.2197E-4 \\ \hline
\textbf{${p_N}$}		& 1.185647 	& 0.828102 	& 0.880824 	& 0.367916 	& 1.404220 	& 0.852152 	& 1.195077 	&  - 		\\ \hline
\hline
\multicolumn{9}{|c|}{Iterative method} \\ \hline
			$h$				& $1/32$   & $1/64$   & $1/128$  & $1/256$  & $1/512$  & $1/1024$ & $1/2048$ & $1/4096$ \\ \hline
\textbf{$\varepsilon$} 	& 0.397955 & 0.228031 & 0.119002 & 0.061311 & 0.033202 & 0,017016 & 0,008611 & 0,004350 \\ \hline
\textbf{${D_N}$} 		& 0.169924 & 0.109029 & 0.057691 & 0.031513 & 0.016186 & 0.008405 & 0.004261 & 0.002151 \\ \hline
\textbf{${p_N}$}		& 0.640177 & 0.918293 & 0.872399 & 0.961200 & 0.945426 & 0.980055 & 0.986184 &  - 		\\ \hline
\end{tabular}
\end{center}
\label{ms_tab3}
\end{figure}

\subsection{Nonlinear equation}
Let us consider the following equation with known solution $\bar{x}(t)=t+\pi$
$$\int\limits_0^{t/8} (t-s) \sin x(s)\;ds+\int\limits_{t/8}^{t/4} t\;(2 \cos x(s))\;ds +\int\limits_{t/4}^{t} (-1)(\sin^2 x(s) + 1)\;ds=$$
$$=-\frac{17t}{8} + \frac{7t}{8}\cos\left(\frac{t}{8}\right) + (1+2t)\sin\left(\frac{t}{8}\right) -2 t \sin\left(\frac{t}{4}\right) - \frac{1}{4}\sin\left(\frac{t}{2}\right) + \frac{1}{4}\sin(2t),\; t\in [0,\,2].$$
\begin{figure}[htb!] %T2
%{Computed maximum pointwise errors $\varepsilon^N,$  two-mesh difference $D^N$ and order of convergence $p^N$ in double precision arithmetic applied to problem 1 for various values of $h=1/N$}
\begin{flushright}
%{Tab. 4.1 Estimations for the 4th example}\vspace{-2mm}
{Tab. 5.2.1 The errors for various stepsizes $h$ of the 4st example.}\vspace{-2mm}
\end{flushright}
\begin{center}
\begin{tabular}{|c|c|c|c|c|c|c|c|c|}	
\hline
\multicolumn{9}{|c|}{Piecewise constant approximation} \\ \hline
		$h$				& $1/32$   & $1/64$   & $1/128$  & $1/256$  & $1/512$  & $1/1024$ & $1/2048$ & $1/4096$ \\ \hline
\textbf{$\varepsilon$} 	& 0.596074 & 0.301323 & 0.130566 & 0.065174 & 0.029381 & 0.014691 & 0.007345 & 0.003672 \\ \hline
\end{tabular}
\end{center}
\label{ms_tab1}
\end{figure}
Tab. \ref{ms_tab1}.1 shows computed maximum pointwise errors $\varepsilon^N=\max\limits_{0\leq i \leq N}|\bar{x}(t_i)-x^h(t_i)|$ for various $h$.

%---------------------------------------------------------------------
%----------------------  S E C T I O N 5 -----------------------------
%---------------------------------------------------------------------
\section{Conclusion}

In this article we proposed the numerical method for solution of the novel class of weakly regular linear and nonlinear Volterra integral
equations of the first kind. We outlined the main
results for this class of equation derived in our previous works. The main contribution
of this paper are a generic numerical methods designed for solution of such weakly regular equations.
The direct numerical methods employe the midpoint quadrature rule and have the
the $\mathcal{O}(1/N)$ and $\mathcal{O}(1/N^2)$ orders of accuracy.  The illustrative examples demonstrate
the efficiency of proposed methods.



%% References


\bibliographystyle{elsarticle-num}


 
%%%%%%%%%%%%%%%%%%%%%%%%

\end{document}